\documentclass[amstex,16pt, amssymb]{article}

\usepackage{mathtext}
\usepackage[cp1251]{inputenc}
\usepackage[T2A]{fontenc}
\usepackage[dvips]{graphicx}
\usepackage{amsmath}
\usepackage{amssymb}
\usepackage{amsxtra}
\usepackage{latexsym}
\usepackage{ifthen}
\usepackage{amsthm}

\numberwithin{equation}{section}

\begin{document}

\title{Hilbert and Poincare problems for semi-linear \\ equations in domains with rectifiable boundaries }

\author{Vladimir Ryazanov}

\date{}

\maketitle

\begin{abstract}

Recall that the research of boundary-value problems with arbitrary
measurable data is due to the famous dissertation of Luzin where he
has studied the Dirichlet problem for harmonic functions in the unit
disk.

In the last paper \cite{R7}, it was studied Hilbert, Poincare and
Neumann boundary-value problems with arbitrary measurable data for
generalized analytic functions and generalized harmonic functions
with applications to the relevant problems of mathematical physics.

The present paper is devoted to the study of the boundary-value
problems with arbitrary measurable boundary data in domains with
rectifiable boun\-da\-ries for the corresponding semi-linear
equations with suitable nonlinear sources.

For this purpose, here it is constructed completely continuous
operators generating nonclassical solutions of the Hilbert and
Poincare boundary-value problems with arbitrary measurable data for
the Vekua type equations and the Poisson equations, respectively.

On this base, it is first proved the existence of solutions of the
Hilbert boundary-value problem with arbitrary mea\-su\-rab\-le data
in any domains with rectifiable boundaries for the nonlinear
equations of the Vekua type.

It is necessary to note that our approach is based on the geometric
interpretation of boundary values as angular (along nontangential
paths) limits in comparison with the classical variational approach
in PDE.

The latter makes it is also possible to obtain the theorem on the
exi\-sten\-ce of nonclassical solutions of the Poincare
boundary-value problem on the directional derivatives and, in
particular, of the Neumann problem with arbitrary mea\-su\-rab\-le
data to the Poisson equations with nonlinear sources in Jordan
domains with rectifiable boundaries.

As consequences, then it is given a series of applications of these
results to some problems of mathematical physics describing such
phenomena as diffusion with physical and chemical absorption, plasma
states and stationary burning.

\end{abstract}




\normalsize \baselineskip=18.5pt

\vskip 1cm


\section{Introduction}
Recall that the boundary-value problems for analytical functions and
their ge\-ne\-ra\-li\-za\-tions go to the famous Riemann
dissertation (1851) and to the following works of Hilbert (1904,
1912, 1924) and Poincare (1910).

The research of boundary-value problems with arbitrary measurable
data is due to the known dissertation of Luzin where he has studied
the corresponding Dirichlet problem for harmonic functions in the
unit disk $\mathbb D := \{ z\in\mathbb C: |z|<1\}$.

In this connection, recall that the following deep result of Luzin
was one of the main theorems of his dissertation, see e.g. his paper
\cite{L1}, dissertation \cite{L2}, p. 35, and its reprint \cite{L},
p. 78, adopted to the segment $[0,2\pi]$.

\medskip

{\bf Theorem A.} {\it\, For any measurable function $\varphi :
[0,2\pi]\to \mathbb R$, there is a continuous function $\Phi :
[0,2\pi]\to \mathbb R$ such that $\Phi^{\prime}=\varphi$ a.e. on
$[0,2\pi]$.}

\medskip

Just on the basis of Theorem A, Luzin proved the next significant
result of his dissertation, see e.g. \cite{L}, p. 87.

\medskip

{\bf Theorem B.} {\it\, Let $\varphi:\mathbb R\to\mathbb R$ be a
$2\pi -$periodic measurable function. Then there is a harmonic
function $u$ in $\mathbb D$ such that $u(z)\to \varphi(\vartheta)$
for a.e. $\vartheta\in\mathbb R$ as $z\to e^{i\vartheta}$ along any
nontangential path.}

\medskip

Note that the Luzin dissertation was later on published only in
Russian as book \cite{L} with comments of his pupils Bari and
Men'shov already after his death. A part of its results was also
printed in Italian \cite{Lu}. However, Theorem A was published in
English in the Saks book \cite{S} as  Theorem VII(2.3). Hence
Frederick Gehring in \cite{Ge} has rediscovered Theorem B and his
proof on the basis of Theorem A in fact coincided with the original
proof of Luzin.

\medskip

Corollary 5.1 in \cite{R1} has strengthened Theorem B as the next,
see also \cite{R2}.

\medskip

{\bf Theorem C.} {\it\, For each measurable function $\varphi
:\partial\mathbb D\to \mathbb R$, the space of all harmonic
functions $u:\mathbb D\to\mathbb R$ with the angular (along
nontangential paths) limits $\varphi(\zeta)$ for a.e.
$\zeta\in\partial\mathbb D$ has the infinite dimension.}

\medskip

The latter was key to establish the following result on the
existence of nonclassical solutions to the Hilbert boundary-value
problem for analytic functions in \cite{R1}, Theorems 3.1 and Remark
5.2.

\medskip

{\bf Theorem D.} {\it\, Let $D$ be a Jordan domain in $\mathbb C$
with a rectifiable boundary and let $\lambda:\partial D\to\mathbb
C$, $|\lambda (\zeta)|\equiv 1$, and $\varphi:\partial D\to\mathbb
R$ be measurable functions. Then the space of analytic functions $f:
D\to\mathbb C$ with the angular limits
\begin{equation}\label{eqLIM} \lim\limits_{z\to\zeta}\ \mathrm {Re}\
\{\overline{\lambda(\zeta)}\cdot f(z)\}\ =\ \varphi(\zeta)
\quad\quad\quad\mbox{for}\ \ \mbox{a.e.}\ \ \ \zeta\in\partial D
\end{equation}
has the infinite dimension.}

\medskip

In turn, on the base of the latter result, it was derived the
corresponding theo\-rems  on the existence of nonclassical solutions
to the Poincare boundary-value problem on the directional
derivatives and, in particular, to the Neumann problem for harmonic
functions in \cite{R3}, Theorems 3, 4 and 5.

\medskip

{\bf Theorem E.} {\it\, Let $D$ be a Jordan domain in $\mathbb C$
with a rectifiable boundary, $\nu:\partial D\to\mathbb C$, $|\nu
(\zeta)|\equiv 1$, and $\varphi:\partial D\to\mathbb R$ be
measurable functions. Then the space of harmonic functions $u:
D\to\mathbb R$ with the angular limits
\begin{equation}\label{eqLIMIT-R} \lim\limits_{z\to\zeta}\ \frac{\partial u}{\partial \nu}\ (z)\ =\
\varphi(\zeta)\quad\quad\quad\mbox{for}\ \ \mbox{a.e.}\ \ \
\zeta\in\partial D\end{equation} has the infinite dimension.}

\medskip

Here we apply the standard designation for the directional
derivative
\begin{equation}\label{eqPOSITIVE}
\frac{\partial u}{\partial \nu}\ (z)\ :=\ \lim_{t\to 0}\
\frac{u(z+t\cdot\nu)-u(z)}{t}\ .
\end{equation}

{\bf Theorem F.} {\it\, Let $D$ be a Jordan domain in $\mathbb C$
with a rectifiable boundary, $\varphi:\partial D\to\mathbb R$ be a
measurable function and let $n=n(\zeta)$ denote the unit interior
normal to $\partial D$ at a point $\zeta$. Then the space of
harmonic functions $u: D\to\mathbb R$ such that, for a.e.
$\zeta\in\partial D$, there exist

\medskip

1) a finite limit along the normal $n(\zeta)$
\begin{equation}\label{eqLIMIT1-R}
u(\zeta)\ :=\ \lim\limits_{z\to\zeta}\ u(z)\ ,\end{equation}

2) the normal derivative
\begin{equation}\label{eqNORMAL-R}
\frac{\partial u}{\partial n}\ (\zeta)\ :=\ \lim_{t\to +0}\
\frac{u(\zeta+t\cdot n)-u(\zeta)}{t}\ =\ \varphi(\zeta)\ ,
\end{equation}

3) the angular limit (along nontangential paths)
\begin{equation}\label{eqLIMIT2-R} \lim\limits_{z\to\zeta}\ \frac{\partial u}{\partial n}\ (z)\ =\
\frac{\partial u}{\partial n}\ (\zeta)\ ,\end{equation} has the
infinite dimension.}

\medskip

In the last paper \cite{R7}, it was studied Hilbert, Poincare and
Neumann boundary-value problems with arbitrary measurable data for
the so--called generalized analytic functions and generalized
harmonic functions with sources and given applications to relevant
problems of mathematical physics.

In this connection, let us recall that the monograph \cite{Ve} was
devoted to the {\bf generalized analytic functions}, i.e.,
continuous complex valued functions $h(z)$ of the complex variable
$z=x+iy$ in $W^{1,1}_{\rm loc}$ satisfying equations of the form
\begin{equation}\label{eqG}
\partial_{\bar z}h\ +\ ah\ +\ bh\ =\ c\ ,\ \ \ \ \ \partial_{\bar
z}\ :=\ \frac{1}{2}\left(\ \frac{\partial}{\partial x}\ +\
i\cdot\frac{\partial}{\partial y}\ \right)\ ,
\end{equation}
where it was assumed that the complex valued functions $a,b$ and $c$
belong to the class $L^{p}$ with some $p>2$ in the corresponding
domain $D\subseteq \mathbb C$.


The paper \cite{R7} contained Theorem 1 on the existence of
nonclassical solutions of the Hilbert boundary-value problem with
arbitrary measurable boundary data for {\bf generalized analytic
functions with sources} $\bf g$, when $a\equiv 0\equiv b$,
\begin{equation}\label{eqS}
\partial_{\bar z}h(z)\ =\ g(z)
\end{equation}
with the real valued functions $g$ in the class $L^{p}$, $p>2$.

Moreover, the paper \cite{R7} included Theorem 6 (Corollary 6)  on
the existence of continuous solutions in $W^{2,p}_{\rm loc}$ to the
Poincare (Neumann) boundary-value problem with arbitrary measurable
boundary data for {\bf generalized harmonic functions with sources}
$\bf G$ in $L^p$, $p>2$, satisfying the Poisson equations
\begin{equation}\label{eqSSS}
\triangle U(z)\ =\ G(z)\ .
\end{equation}

The present paper is devoted to the study of the boundary-value
problems with arbitrary measurable boundary data in domains with
rectifiable boun\-da\-ries for the corresponding semi-linear
equations with suitable nonlinear sources.

Namely, the first part of the paper is devoted to the Hilbert
boundary-value problem with arbitrary measurable boundary data in
Jordan domains $D$ with rectifiable boundaries for the nonlinear
Vekua type equations of the form
\begin{equation} \label{eqQUASILINEARWintr}
\partial_{\bar z}f(z)\ =\ h(z)\cdot q(f(z))\ \ \ \ \ \ \ \mbox{a.e. in $D$}\ ,
\end{equation}
where $h: D\to\mathbb C$ is a function in the class $L^p(D)$ for
$p>2$ and $q:\mathbb C\to\mathbb C$ is a continuous function with
\begin{equation} \label{eqAPRIORYW}
\lim\limits_{w\to \infty}\ \frac{q(w)}{w}\ =\ 0\ .
\end{equation}

The second part of the paper is devoted to the Poincare (and
Neumann) boundary-value problem with arbitrary measurable boundary
data in Jordan domains $D$ with rectifiable boundaries for the
nonlinear Poisson equations \begin{equation} \label{eqQUASILINEARWH}
\triangle U(z)\ =\ H(z)\cdot Q(U(z))\ \ \ \ \ \ \ \mbox{a.e. in
$D$}\ ,
\end{equation} where
$H: D\to\mathbb R$ is a function in the class $L^p(D)$ for $p>2$ and
$Q:\mathbb R\to\mathbb R$ is a continuous function with
\begin{equation} \label{eqAPRIORYWHintr}
\lim\limits_{t\to \infty}\ \frac{Q(t)}{t}\ =\ 0\ .
\end{equation}

For this purpose, it is first established the existence of
completely con\-ti\-nuo\-us operators generating nonclassical
solutions of the Hilbert and Poincare boundary-value problems with
arbitrary measurable data for the equations of the Vekua and Poisson
types (\ref{eqS}) and (\ref{eqSSS}), respectively.

Finally, the third part includes a series of applications of the
results on the Poincare and Neumann boundary-value problems to some
nonlinear equations of mathematical physics modeling, for instance,
such phenomena as physical and che\-mi\-cal absorption with
diffusion, plasma states, stationary burning etc.

Here we use the geometric interpretation of boundary values as
angular (along nontangential paths) limits that is a traditional
tool in the geometric function theory, see e.g. monographs
\cite{Du}, \cite{Ko}, \cite{L}, \cite{Po} and \cite{P}.

\section{Definitions and preliminary remarks}

First of all, recall that a {\bf completely continuous} mapping from
a metric space $M_1$ into a metric space $M_2$ is defined as a
continuous mapping on $M_1$ which takes bounded subsets of $M_1$
into relatively compact ones of $M_2$, i.e. with compact closures in
$M_2$. When a continuous mapping takes $M_1$ into a relatively
compact subset of $M_1$, it is nowadays said to be {\bf compact} on
$M_1$.

Note that the notion of completely continuous (compact) operators is
due essentially, in the special space of bilinear forms in $l_2$, to
Hilbert who requires the operator to map weakly convergent sequences
into strongly convergent sequences that, in reflexive spaces, is
equivalent to Definition VI.5.1 for the Banach spaces in \cite{DS}
which is due to F. Riesz, see also the comments of Section VI.12 in
\cite{DS}. The latter just coincides with the above definition in
the special case.

Recall more some definitions and the fundamental result of the
celebrated paper \cite{LS}. Leray and Schauder extend as follows the
Brouwer degree to compact perturbations of the identity $I$ in a
Banach space $B$, i.e. a complete normed linear space. Namely, given
an open bounded set $\Omega\subset B$, a compact mapping $F: B\to B$
and $z \notin \Phi(\partial \Omega)$, $\Phi :=I-F$, the {\bf
(Leray–Schauder) topological degree} $\deg\, [\Phi,\Omega, z]$ of
$\Phi$ in $\Omega$ over $z$ is constructed from the Brouwer degree
by approximating the mapping $F$ over $\Omega$ by mappings
$F_{\varepsilon}$ with range in a finite-dimensional subspace
$B_{\varepsilon}$ (containing $z$) of $B$. It is showing that the
Brouwer degrees $\deg\, [\Phi_{\varepsilon} ,\Omega_{\varepsilon},
z]$ of $\Phi_{\varepsilon}:=I_{\varepsilon} - F_{\varepsilon}$,
$I_{\varepsilon}:=I|_{B_{\varepsilon}}$, in
$\Omega_{\varepsilon}:=\Omega\cap B_{\varepsilon}$ over $z$
stabilize for sufficiently small positive $\varepsilon$ to a common
value defining $\deg\, [\Phi,\Omega, z]$ of $\Phi$ in $\Omega$ over
$z$.


This topological degree “algebraically counts” the number of fixed
points of $F(\cdot)-z$ in $\Omega$ and conserves the basic
properties of the Brouwer degree as ad\-di\-ti\-vi\-ty and homotopy
invariance. Now, let $a$ be an isolated fixed point of $F$. Then the
{\bf local (Leray–Schauder) index} of $a$ is defined by ${\rm ind}\,
[\Phi, a] := \deg [\Phi,B(a, r), 0]$ for small enough $r > 0$. ${\rm
ind}\, [\Phi, 0]$ is called by {\bf index} of $F$. In particular, if
$F\equiv 0$, correspondingly, $\Phi\equiv I$, then the index of $F$
is equal to $1$.


Let us formulate the main result in \cite{LS}, Theorem 1, see also
the survey \cite{Ma}.

\bigskip

{\bf Proposition 1.} {\it Let $B$ be a Banach space, and let
$F(\cdot,\tau):B\to B$ be a family of operators with $\tau\in[0,1]$.
Suppose that the following hypotheses hold:

\medskip

{\rm {\bf (H1)}} $F(\cdot,\tau)$ is completely continuous on $B$ for
each $\tau\in[0,1]$ and uniformly continuous with respect to the
parameter $\tau\in[0,1]$ on each bounded set in $B$;

{\rm {\bf (H2)}} the operator $F:=F(\cdot,0)$ has finite collection
of fixed points whose total index is not equal to zero;

{\rm {\bf (H3)}} the collection of all fixed points of the operators
$F(\cdot,\tau)$, $\tau\in[0,1]$,  is bounded in $B$.

Then the collection of all fixed points of the family of operators
$F(\cdot,\tau)$ contains a continuum along which $\tau$ takes all
values in  $[0,1]$.}

\medskip

Let us go back to the discussion of the results of Luzin in
Introduction.

\medskip

{\bf Remark 1.} Applying the Cantor ladder type functions, namely,
continuous nondecreasing functions $C:[0,2\pi]\to\mathbb R$ with
$C(0)=0$, $C(2\pi)=1$ and $C^{\prime}(t)=0$ for a.e. $t\in[0,2\pi]$,
see e.g. Section 8.15 in \cite{GO}, we may assume in Theorem A that
$\Phi(0)=0=\Phi(2\pi)$. On the same base, using uniform continuity
of the function $\Phi$ on $[0,2\pi]$ and applying sequentially
fragmentations of the segment to arbitrarily small parts, we may
assume in Theorem A that $|\Phi(t)|<\varepsilon$ for every
prescribed $\varepsilon>0$ and, in particular, that $|\Phi(t)|<1$
for all $t\in[0,2\pi]$. Thus, in view of arbitrariness of
$\varepsilon>0$, there is the infinite collection of such $\Phi$ for
each $\varphi$. Furthermore, applying series of pair of
(nondecreasing and nonincreasing) functions of the Cantor ladder
type on the segments $[2^{-(k+1)}\pi,2^{-k}\pi]$, $k=1,2,\ldots$ it
is easy to see that the space of such functions $\Phi$ has the
infinite dimension.

\medskip

By the proof of Theorem B, see \cite{L2}, \cite{L} or \cite{Ge},
$u(z)\ =\ \frac{\partial}{\partial\vartheta}\ U(z)$, where
\begin{equation}\label{Poisson}
U(re^{i\vartheta})\ =\ \frac{1}{2\pi}\
\int\limits_{0}\limits^{2\pi}\frac{1-r^2}{1-2r\cos(\vartheta-t)+r^2}\
\Phi(e^{it})\ dt\ ,
\end{equation}
i.e., for a function $\Phi$ from Theorem A, $u$ can be calculated in
the explicit form
\begin{equation}\label{eqG}
u(re^{i\vartheta})\ =\ -\ \frac{r}{\pi}\
\int\limits_{0}\limits^{2\pi}\frac{(1-r^2)\sin(\vartheta-t)}{(1-2r\cos(\vartheta-t)+r^2)^2}\
\Phi(e^{it})\ dt\ . \end{equation}

{\bf Remark 2.} Later on, it was shown by Theorems 3 in \cite{R4}
that the Luzin harmonic functions $u(z)$ can be represented as the
{\bf Poisson--Stieltjes integrals}
\begin{equation}\label{eqPS1}
\mathbb U_{\Phi}(z)\ =\ \frac{1}{2\pi}\
\int\limits_{-\pi}\limits^{\pi} P_r(\vartheta -t)\ d\,\Phi(e^{it})\
\ \ \ \forall\ z=re^{i\vartheta}, \ r\in(0,1)\ ,\
\vartheta\in[-\pi,\pi]\ ,
\end{equation}
where $P_r(\Theta) = (1-r^2)/(1-2r\cos\Theta +r^2), r<1,
\Theta\in\mathbb R,$ is the {\bf Poisson kernel}.

The corresponding analytic functions in $\mathbb D$ with the real
parts $u(z)$ can be represented  as the corresponding {\bf
Schwartz--Stieltjes integrals}
\begin{equation}\label{eqS1}
\mathbb S_{\Phi}(z)\ =\ \frac{1}{2\pi}\ \int\limits_{\partial\mathbb
D}\frac{\zeta +z}{\zeta -z}\ d\,\Phi(\zeta)\ ,\ \ \ \ z\in \mathbb
D\ ,
\end{equation}
because of the Poisson kernel is the real part of the (analytic in
the variable $z$) {\bf Schwartz kernel} $(\zeta +z)/(\zeta -z)$.
Integrating (\ref{eqS1}) by parts, see Lemma 1 and Remark 1 in
\cite{R4}, we obtain also the more convenient form of the
representation
\begin{equation}\label{eqS2}
\mathbb S_{\Phi}(z)\ =\ \frac{z}{\pi}\ \int\limits_{\partial\mathbb
D}\frac{\Phi(\zeta)}{(\zeta -z)^2}\ d\,\zeta\ ,\ \ \ \ z\in \mathbb
D\ .
\end{equation}

\section{On completely continuous Hilbert operators}

Recall that in paper \cite{R7}, we considered {generalized analytic
functions $f$ with sources $g\in L^p, p>2,$} in the class
$W^{1,1}_{\rm loc}$ that satisfy the equation
\begin{equation}\label{eqV}  \frac{\partial f}{\partial{\bar
z}}\ =\ g\ ,\ \ \ \ \ \frac{\partial}{\partial{\bar z}}\ :=\
\frac{1}{2}\left(\ \frac{\partial}{\partial x}\ +\
i\cdot\frac{\partial}{\partial y}\ \right)\ ,\ \ \ z=x+iy\ ,
\end{equation}
and studied for them the Hilbert boundary-value problem in Jordan
domains with rectifiable boundaries under arbitrary boun\-dary data
that are measurable over the natural parameter.

In particular, Theorem 1 in \cite{R7} stated that, for arbitrary
measurable functions $\lambda:\partial\mathbb D\to\mathbb{C},\:
|\lambda(\zeta)|\equiv1$, and $\varphi:\partial\mathbb
D\to\mathbb{R}$, there exist generalized analytic functions $f:
D\to\mathbb C$  with any source $g:\mathbb D\to \mathbb R$ in the
class $L^p(\mathbb D)$, $p>2$, that have the angular limits
\begin{equation}\label{eqLIMH} \lim\limits_{z\to\zeta}\ \mathrm
{Re}\ \left\{\, \overline{\lambda(\zeta)}\cdot f(z)\, \right\}\ =\
\varphi(\zeta) \quad\quad\quad \mbox{a.e.\ on $\ \partial\mathbb
D$}\ .\end{equation} Furthermore, the space of such functions $f$
has the infinite dimension.

Thus, the Hilbert boundary-value problem always has many solutions
in the given sense for each such coefficient $\lambda$, boundary
date $\varphi$ and source $g$. Of course, axiom of choice by Zermelo
makes it possible to choose one of such correspondence named further
as a Hilbert operator but the latter with such a random choice can
be completely discontinuous. Later on, to apply the approach of
Leray-Schauder for extending Theorem 1 in \cite{R7} to the
generalized analytic functions, satisfying nonlinear equations of
the Vekua type, we need just the complete continuity of such
correspondence.


So, let us construct a completely continuous Hilbert operator
generating generalized analytic functions with sources $g:\mathbb
D\to\mathbb C$ in the class $L^p$, $p>2$, and the boundary condition
(\ref{eqLIMH}) for prescribed measurable functions $\lambda$ and
$\varphi$.

For this purpose, let us first consider the known linear singular
operator
\begin{equation}\label{eqRELATION} T_g(z)\ :=\ \frac{1}{\pi}\int\limits_{\mathbb C} g(w)\ \frac{d\,
m(w)}{z-w}\ , \end{equation} where we assume that $g$ is extended by
zero outside of $\mathbb D$.

\medskip

{\bf Remark 3.} By Theorem 1.14 in \cite{Ve} the function $T_g$ has
the generalized derivative by Sobolev $\partial T_g/\partial\bar
z=g$ if $g\in L^1(\mathbb D)$. Moreover, by Theorem 1.36 in
\cite{Ve} the function $T_g\in W^{1,p}_{\rm loc}$ if $g\in
L^p(\mathbb D)$, $p>1$.

Furthermore, if $g\in L^p(\mathbb D)$, $p>2$, then by Theorem 1.19
in \cite{Ve}
\begin{equation}\label{eqH1}
|T_g(z)|\ \leq\ M_1\,\| g\|_p\ \ \ \ \ \forall\ z\in\mathbb C\ ,
\end{equation}
\begin{equation}\label{eqH1}
|T_g(z_1) - T_g(z_2)|\ \leq\ M_2\,\| g\|_p\, |z_1-z_2|^{\alpha}\ \ \
\ \ \forall\ z_1,z_2\in\mathbb C\ ,
\end{equation}
where the constants $M_1$ and $M_2$ depend only on $p>2$, and
$\alpha=(p-2)/p$. Thus, the linear operator $T_g$ is completely
continuous on compact sets in $\mathbb C$ and, in particular, on
$\overline{\mathbb D}$ by Arzela-Ascoli theorem, see e.g. Theorem
IV.6.7 in \cite{DS}.

\bigskip

Next, since $T_g$ is continuous, we have the measurable boundary
function
\begin{equation}\label{eqH} \varphi_g(\zeta)\ :=\ \lim\limits_{z\to\zeta}\ \mathrm
{Re}\ \left\{\, \overline{\lambda(\zeta)}\cdot T_g(z)\, \right\}\ =\
\mathrm {Re}\ \left\{\, \overline{\lambda(\zeta)}\cdot T_g(\zeta)\,
\right\}\ ,\ \ \ \ \forall\,\zeta\in\partial\mathbb D\
.\end{equation}

Thus, the generalized analytic functions $f$ with the source $g$
satisfying the Hilbert condition (\ref{eqLIMH}) can be get as the
sums $f=T_g+{\cal C}$ with analytic functions $\cal C$ satisfying,
in the sense of angular limits, the Hilbert boundary condition
\begin{equation}\label{eqLIMDD} \lim\limits_{z\to\zeta}\ \mathrm
{Re}\ \{\overline{\lambda(\zeta)}\cdot {\cal C}(z)\}\ =\
\psi(\zeta):=\varphi(\zeta)-\varphi_g(\zeta) \quad\quad\quad
\mbox{a.e.\ on\ $\partial\mathbb D$}\ .\end{equation}

In turn, by the construction of Theorem 2.1 in \cite{R1}, such
analytic functions ${\cal C}$ can be obtained as the products of 2
analytic functions ${\cal A}$ and ${\cal B}$. The first
\begin{equation}\label{eqA} {\cal A}(z)\ =\ e^{ia(z)}\ , \ \ \ a(z)\
:=\ \frac{1}{2\pi i}\ \int\limits_{\partial\mathbb
D}\alpha_{\lambda}(\zeta)\ \frac{z+\zeta}{z-\zeta}\
 \frac{d\zeta}{\zeta}\ , \ \ \ \ \ z\in\mathbb D\ ,
\end{equation} with $\alpha_{\lambda}(\zeta):=\arg\lambda(\zeta)$,
where $\arg \omega$  is the principal value of the argument of
$\omega\in\mathbb C$, $|\omega|=1$, i.e., the number $\alpha\in
(-\pi,\pi]$ such that $\omega=e^{i\alpha}$; the second one
\begin{equation}\label{eqS3} {\cal B}(z)\ =\ \mathbb S_{\Psi}(z)\ =\
\frac{z}{\pi}\ \int\limits_{\partial\mathbb
D}\frac{\Psi(\zeta)}{(\zeta -z)^2}\ d\,\zeta\ ,\ \ \ \ z\in \mathbb
D\ ,
\end{equation} see Remark 2, where $\Psi$ is an antiderivative of the function
$\psi e^{\beta}$ in Theorem A, $\beta(\zeta)$ is the angular limit
of ${\mathrm {Im}}\ a(z)$ as $z\to\zeta\in\partial\mathbb D$ a.e.,
see e.g. Corollary 4 in \cite{R4}, ${\mathrm {Re}}\
a(z)\to\alpha_{\lambda}(\zeta)$ as $z\to\zeta\in\partial\mathbb D$
a.e., see e.g. Corollary IX.1.1 in \cite{Gol}.

Thus, analytic functions $\cal C$ can be represented in more
convenient form
\begin{equation}\label{eqS4}
{\cal C}(z)\ =\ {\cal A}(z)\cdot[\,\mathbb S_{\Phi}(z)\, -\, \mathbb
S_{\Phi_g}(z)\,]\ ,
\end{equation}
where $\Phi$ and $\Phi_g$ are antiderivatives of $\varphi e^{\beta}$
and $\varphi_g e^{\beta}$ as functions of $\vartheta\in[0,2\pi]$ in
Theorem A, respectively. Note that the analytic functions $\cal A$
and $\mathbb S_{\Phi}$ do not depend on the sources $g$ at all. Let
us choose the function $\Phi_g$ in a suitable way.

From this point on, we demand that all sources $g$ have compact
supports in the unit disk and belong to a disk $\mathbb D_{\rho}:=\{
z\in\mathbb C: |z|\leq\rho\}$ with a radius $\rho \in (0,1)$. Then
the function $T_g(z)$, $z\in\mathbb C$, is analytic in a
neighborhood of the unit circle $\partial\mathbb D$ and, in
particular, $T_g(\zeta)$ is differentiable in $\vartheta\in\mathbb
R$, $\zeta=e^{i\vartheta}$. Moreover, we have that, for all
$\zeta=e^{i\vartheta}$, $\vartheta\in[0,2\pi]$,
\begin{equation}\label{eqRELATION3} \{ T_g\}_{\vartheta}(\zeta)\ =\ i\zeta
T_g^{\prime}(\zeta)\ =\ \frac{\zeta}{\pi i}\int\limits_{\mathbb
D_{\rho}} g(w)\ \frac{d\, m(w)}{(\zeta-w)^2}\ \ \ \ \ \ \ \ \forall\
\zeta\in\partial\mathbb D\ .
\end{equation}
Let us denote by $\Lambda$ an antiderivative for the function
$\overline{\lambda}e^{\beta}$ as a function of
$\vartheta\in[0,2\pi]$ in Theorem A, see also Remark 1.

Then the following function $\Phi_g$ is an antiderivative for the
function $\varphi_ge^{\beta}$:
\begin{equation}\label{eqS5} \Phi_g(\zeta)\ :=\
\mathrm {Re}\ \left\{\, \Lambda(\zeta) T_g(\zeta) -
\int\limits_{0}^{\vartheta} \Lambda(\xi)\{ T_g\}_{\theta}(\xi)\,
d\,\theta + S(\vartheta) \, \right\}\ ,\end{equation} where
$S:[0,2\pi]\to\mathbb C$ is either zero or a singular function of
the form
\begin{equation}\label{eqS6} S(\vartheta)\ :=\ C(\vartheta)\int\limits_{0}^{2\pi}
\Lambda(\xi)\{ T_g\}_{\theta}(\xi )\, d\,\theta \ , \ \ \ \
\zeta=e^{i\vartheta},\ \xi=e^{i\theta} ,\ \vartheta,\
\theta\in[0,2\pi]\ ,\end{equation} with a singular function
$C:[0,2\pi]\to[0,1]$ of the Cantor ladder type, see Section 8.15 in
\cite{GO}, i.e., $C$ is continuous, nondecreasing, $C(0)=0$,
$C(2\pi)=1$ and $C^{\prime}=0$ a.e. on $[0,2\pi]$.

Let us show that the Hilbert operator ${\cal H}^*_g$ generated by
the sums $T_g+{\cal C}$ under the given choice of $\Phi$ and
$\Phi_g$ in (\ref{eqS4}) is completely continuous on compact sets in
$\mathbb D$. Recall that the analytic functions $\cal A$ and
$\mathbb S_{\Phi}$ in the representation (\ref{eqS4}) of ${\cal C}$
do not depend on the sources $g$. Hence by Remark 3, it remains to
show that the linear operator $\mathbb S_{\Phi_g}$ is completely
continuous.

Indeed, by the construction of $\Phi_g$ in (\ref{eqS5}) and
relations (\ref{eqRELATION}) and (\ref{eqRELATION3})
\begin{equation}\label{eqEST}
|\Phi_g(\zeta)|\ \leq\ \frac{1}{\pi}\cdot\frac{\|g\|_1}{1-\rho}\ +\
2\cdot\frac{\|g\|_1}{(1-\rho)^2}\ \leq\ c_{\rho}\cdot\|g\|_1\ \leq\
C_{\rho}\cdot\|g\|_p\ \ \ \forall\ \zeta\in\partial\mathbb D
\end{equation}
with $c_{\rho}=3/(1-\rho)^2$ and  $C_{\rho}=3\pi/(1-\rho)^2$,
respectively. Hence, by (\ref{eqS2})
\begin{equation}\label{eqEST2}
|\mathbb S_{\Phi_g}(z)|\ \leq\ C_{\rho,r}\cdot\|g\|_p\ ,\ \ \
\forall\ z\in\mathbb D_r\ ,\ r\in(0,1)\ ,
\end{equation}
\begin{equation}\label{eqEST3}
|\mathbb S_{\Phi_g}(z_1)-\mathbb S_{\Phi_g}(z_2)|\ \leq\
C^*_{\rho,r}\cdot\|g\|_p\cdot |z_1-z_2|\ ,\ \ \ \forall\
z_1,z_2\in\mathbb D_r\ ,\ r\in(0,1)\ ,
\end{equation}
where the constants $C_{\rho,r}$ and $C^*_{\rho,r}$ depend only on
the radii $\rho$ and $r\in(0,1)$. Thus, the operator $\mathbb
S_{\Phi_g}$ is completely continuous on compact sets in $\mathbb D$
again by the Arzela-Ascoli theorem. Combining it with Remark 3, we
obtain the following conclusion.

\medskip

{\bf Lemma 1.}{\it\, Let $\lambda:\partial\mathbb D\to\mathbb{C},\:
|\lambda(\zeta)|\equiv1$, and $\varphi:\partial\mathbb
D\to\mathbb{R}$ be measurable. Then there is a Hilbert operator
${\cal H}^*_g$ over $g:\mathbb D\to \mathbb C$ in $L^p(\mathbb D)$,
$p>2$, with compact supports in $\mathbb D$, generating generalized
analytic functions $f: \mathbb D\to\mathbb C$ with the sources $g$
and the angular limits
\begin{equation}\label{eqLIMH1} \lim\limits_{z\to\zeta}\ \mathrm
{Re}\ \left\{\, \overline{\lambda(\zeta)}\cdot f(z)\, \right\}\ =\
\varphi(\zeta) \quad\quad\quad \mbox{a.e.\ on $\ \partial\mathbb
D$}\ ,\end{equation} whose restriction to sources $g$ with
$\rm{supp}\,g\subseteq\mathbb D_{\rho}$ is completely continuous
over $\mathbb D_r$ for each $\rho$ and $r\in(0,1).$ }

\medskip

{\bf Remark 4.} Note that the nonlinear operator ${\cal H}^*_g$
constructed above is not bounded except the trivial case $\Phi\equiv
0$ because then ${\cal H}^*_0={\cal A}\cdot\mathbb S_{\Phi}\neq 0$.
However, the restriction of the operator ${\cal H}^*_g$ to $\mathbb
D_r$ under each $r\in(0,1)$ is bounded at infinity in the sense that
$\max\limits_{z\in\mathbb D_r}|{\cal H}^*_g(z)|\leq M\cdot\|g\|_p$
for some $M>0$ and all $g$ with large enough $\| g\|_p$. Note also
that by Remark 1 we are able always to choose $\Phi$ for any
$\varphi$, including $\varphi\equiv 0$, which is not identically $0$
in the unit disk $\mathbb D$.

\bigskip

\section{On Hilbert problem for semi-linear equations}

In this section we study the solvability of the Hilbert
boundary-value problem for nonlinear equations of the Vekua type
$\partial_{\bar z}f(z) = h(z)q(f(z))$ in the unit disk $\mathbb D$.
The Leray--Schauder approach described in Section 2 allows us to
reduce the problem to the study of the corresponding linear equation
from our last paper \cite{R7} on the basis of Lemma 1 in the
previous section on completely continuous Hilbert operator ${\cal
H}^*_g$ and Remark 4 on its boundedness at infinity.

\bigskip

In the proof of the next theorem, the initial operator
$F(\cdot):=F(\cdot,0)\equiv 0$. Hence $F$ has the only one fixed
point (at the origin) and its index is equal to $1$ and, thus,
hypothesis (H2) in Proposition 1 will be automatically satisfied.


{\bf Theorem 1.} {\it Let $\lambda:\partial\mathbb D\to\mathbb{C},\:
|\lambda(\zeta)|\equiv1$, and $\varphi:\partial\mathbb
D\to\mathbb{R}$ be measurable. Suppose that $h: \mathbb D\to\mathbb
C$ is in $L^p(\mathbb D)$ for $p>2$ with compact support in $\mathbb
D$ and $q:\mathbb C\to\mathbb C$ is a continuous function with
\begin{equation} \label{eqAPRIORYW}
\lim\limits_{w\to \infty}\ \frac{q(w)}{w}\ =\ 0\ .
\end{equation}

Then there is $f:\mathbb D\to\mathbb C$ in the class $W^{1,p}_{\rm
loc}\cap C^{\alpha}_{\rm loc}(\mathbb D)$ with $\alpha=(p-2)/p$,
\begin{equation} \label{eqQUASILINEARW}
\partial_{\bar z}f(z)\ =\ h(z)\cdot q(f(z))\ \ \ \ \ \ \ \mbox{a.e. in $\mathbb
D$}\ ,
\end{equation}
with the angular limits
\begin{equation}\label{eqLIMHh} \lim\limits_{z\to\zeta}\ \mathrm
{Re}\ \left\{\, \overline{\lambda(\zeta)}\cdot f(z)\, \right\}\ =\
\varphi(\zeta) \quad\quad\quad \mbox{a.e.\ on $\ \partial\mathbb
D$}\ .\end{equation} }


\begin{proof} If $\| h\|_p=0$ or $\| q\|_C= 0$, then any analytic
function from Theorem 2.1 in \cite{R7} gives the desired solution of
(\ref{eqQUASILINEARW}). Thus, we may assume that $\| h\|_p\neq 0$
and $\| q\|_C\neq 0$. Set $q_*(t)=\max\limits_{|w|\le t}|q(w)|$,
$t\in\mathbb R^+:=[0,\infty)$. Then the function $q_*:\mathbb
R^+\to\mathbb R^+$ is continuous and nondecreasing and, moreover, by
(\ref{eqAPRIORYW})
\begin{equation} \label{eqAPRIORYWT}
\lim\limits_{t\to \infty}\ \frac{q_*(t)}{t}\ =\ 0\ .
\end{equation}

By Lemma 1 and Remark 4 we obtain the family of operators
$F(g;\tau): L_h^{p}(\mathbb D)\to L_h^{p}(\mathbb D)$, where
$L_h^{p}(\mathbb D)$ consists of functions $g\in L^{p}(\mathbb D)$
with supports in the support of $h$,
\begin{equation} \label{eqFORMULA}
F(g;\tau)\ :=\ \tau h\cdot q({\cal H}^*_g) \ \ \ \ \ \ \ \forall\
\tau\in[0,1]
\end{equation}
which satisfies all groups of hypothesis H1-H3 of Theorem 1 in
\cite{LS}, see Proposition 1. Indeed:

H1). First of all, by Lemma 1 the function $F(g;\tau)\in
L_h^{p}(\mathbb D)$ for all $\tau\in[0,1]$ and  $g\in
L_h^{p}(\mathbb C)$ because the function $q({\cal H}^*_g)$ is
continuous and, furthermore, the operators $F(\cdot ;\tau)$ are
completely continuous for each $\tau\in[0,1]$ and even uniformly
continuous with respect to the parameter $\tau\in[0,1]$.

H2). The index of the operator $F(g;0)$ is obviously equal to $1$.

H3). Let us assume that solutions of the equations $g=F(g;\tau)$ is
not bounded in $L_h^{p}(\mathbb D)$, i.e., there is a sequence of
functions $g_n\in L_h^{p}(\mathbb D)$ with $\|g_n\|_p\to\infty$ as
$n\to\infty$ such that $g_n=F(g_n;\tau_n)$ for some
$\tau_n\in[0,1]$, $n=1,2,\ldots$.

However, then by Remark 4 we have that, for some constant $M>0$,
$$
\| g_n\|_p\ \le\ \| h\|_p\ q_*\left(\, M\, \| g_n\|_p\right)\
$$
and, consequently,
\begin{equation} \label{eqEST3W}
\frac{q_*(\, M\, \| g_n\|_p)}{M\, \| g_n\|_p}\ \ge\ \frac{1}{M\, \|
h\|_p}\ >\ 0
\end{equation}
for all large enough $n$. The latter is impossible by condition
(\ref{eqAPRIORYWT}). The obtained contradiction disproves the above
assumption.

Thus, by Theorem 1 in \cite{LS} there is a function $g\in L_h^p(D)$
with $F(g;1)=g$, and by Lemma 1 the function $f:={\cal H}^*_g$ gives
the desired solution of (\ref{eqQUASILINEARW}).
\end{proof}


In particular, choosing $\lambda\equiv 1$ in Theorem 1 we obtain the
following consequence on the Dirichlet problem for the nonlinear
equations of the Vekua type.

\medskip

{\bf Corollary 1.} {\it Let $\varphi:\partial\mathbb D\to\mathbb{R}$
be a measurable function, $h: \mathbb D\to\mathbb C$ be a function
in the class $L^p(\mathbb D)$ for $p>2$ with compact support in
$\mathbb D$ and let $q:\mathbb C\to\mathbb C$ be a continuous
function with condition (\ref{eqAPRIORYW}).

Then there is a function $f:\mathbb D\to\mathbb C$ in the class
$W^{1,p}_{\rm loc}\cap C^{\alpha}_{\rm loc}(\mathbb D)$ with
$\alpha=(p-2)/p$, satisfying equation (\ref{eqQUASILINEARW}) a.e.
with the angular limits
\begin{equation}\label{eqLIMHhD} \lim\limits_{z\to\zeta}\ \mathrm
{Re}\  f(z)\ =\ \varphi(\zeta) \quad\quad\quad \mbox{a.e.\ on $\
\partial\mathbb D$}\ .\end{equation}
}


{\bf Remark 5.} Moreover, by the proof of Theorem 1 $f$ is a
generalized analytic function with a source $g\in L^p(\mathbb D)$,
$f={\cal H}^*_g$, where ${\cal H}^*_g$ is the Hilbert operator
described in the last section (with the simplest analytic function
${\cal A}\equiv 1$ in the case of Corollary 1), Lemma 1, and the
support of $g$ is in the support of $h$ and the upper bound of $\|
g\|_p$ depends only on $\| h\|_p$ and on the function $q$.

In addition, the source $g:\mathbb D\to\mathbb C$ is a fixed point
of the nonlinear operator $\Omega_g:=h\cdot q({\cal
H}^*_g):L^p_h(\mathbb D)\to L^p_h(\mathbb D)$, where
$L_h^{p}(\mathbb D)$ consists of functions $g$ in $L^{p}(\mathbb D)$
with supports in the support of $h$.

\section{On the Hilbert problem in rectifiable domains}

{\bf Theorem 2.} {\it Let $D$ be  a Jordan domain in $\mathbb C$
with a rectifiable boundary, $\lambda:\partial D\to\mathbb{C},\:
|\lambda(\zeta)|\equiv1$, and $\varphi:\partial D\to\mathbb{R}$ be
measurable over natural parameter.

Suppose that $h: D\to\mathbb C$ is in $L^p(D)$ for $p>2$ with
compact support in $D$ and $q:\mathbb C\to\mathbb C$ is a continuous
function with
\begin{equation} \label{eqAPRIORYW5}
\lim\limits_{w\to \infty}\ \frac{q(w)}{w}\ =\ 0\ .
\end{equation}

Then there is $f:D\to\mathbb C$ in the class $W^{1,p}_{\rm loc}\cap
C^{\alpha}_{\rm loc}(D)$ with $\alpha=(p-2)/p$,
\begin{equation} \label{eqQUASILINEARW5}
\partial_{\bar \xi}f(\xi)\ =\ h(\xi)\cdot q(f(\xi))\ \ \ \ \ \ \ \mbox{a.e. in
$D$}\ ,
\end{equation}
and the angular limits
\begin{equation}\label{eqLIMHh5} \lim\limits_{\xi\to\omega}\ \mathrm
{Re}\ \left\{\, \overline{\lambda(\omega)}\cdot f(\xi)\, \right\}\
=\ \varphi(\omega) \quad\quad\quad \mbox{a.e.\ on $\ \partial D$}\
.\end{equation} }


\begin{proof} Let $c$ be a conformal mapping of $D$ onto $\mathbb D$ that exists
by the Riemann mapping theorem, see e.g. Theorem II.2.1 in
\cite{Gol}. Now, by the Caratheodory theorem, see e.g. Theorem
II.3.4 in \cite{Gol}, $c$ is extended to a homeomorphism $\tilde c$
of $\overline{D}$ onto $\overline{\mathbb D}$. Set $c_*=\tilde
c|_{\partial D}$. If $\partial D$ is rectifiable, then by the
theorem of F. and M. Riesz $\mathrm{length}\ c_*^{-1}(E)=0$ whenever
$E\subset\partial\mathbb D$ with $|E|=0$, see e.g. Theorem II.C.1
and Theorems II.D.2 in \cite{Ko}. Conversely, by the Lavrentiev
theorem $|c_*({\cal E})|=0$ whenever ${\cal E}\subset\partial D$ and
$\mathrm{length}\ {\cal E}=0$, see \cite{Lavr}, see also the point
III.1.5 in \cite{P}.

Hence $c_*$ and $c_*^{-1}$ transform measurable sets into measurable
sets. Indeed, every measurable set is the union of a sigma-compact
set and a set of measure zero, see e.g. Theorem III(6.6) in
\cite{S}, and continuous mappings transform compact sets into
compact sets. Thus, functions $\lambda:\partial D\to\mathbb C$,
$|\lambda|\equiv 1$, and $\varphi:\partial D\to\mathbb R$ are
measurable with respect to the natural parameter on $\partial D$ if
and only if the functions $\tilde\lambda=\lambda\circ
c_*^{-1}:\partial\mathbb D\to\mathbb C$ and
$\tilde\varphi=\varphi\circ c_*^{-1}:\partial\mathbb D\to\mathbb R$
are so.

Now, set $\tilde h=h\circ C\cdot \overline{C^{\prime}}$, where $C$
is the inverse conformal mapping to $c$, $C:=c^{-1}:\mathbb D\to D$.
Then it is clear by the hypotheses of Theorem 2 that $\tilde h$ has
compact support in $\mathbb D$ and belongs to the class $L^p(\mathbb
D)$. Consequently, by Theorem 1 there is $\tilde f:\mathbb
D\to\mathbb C$ in the class $W^{1,p}_{\rm loc}\cap C^{\alpha}_{\rm
loc}(\mathbb D)$ with $\alpha=(p-2)/p$,
\begin{equation} \label{eqQUASILINEARW55}
\partial_{\bar z}\tilde f(z)\ =\ \tilde h(z)\cdot q(\tilde f(z))\ \ \ \ \ \ \ \mbox{a.e. in $\mathbb
D$}\ ,
\end{equation}
and the angular limits
\begin{equation}\label{eqLIMHh55} \lim\limits_{z\to\zeta}\ \mathrm
{Re}\ \left\{\, \overline{\tilde\lambda(\zeta)}\cdot\tilde f(z)\,
\right\}\ =\ \tilde\varphi(\zeta) \quad\quad\quad \mbox{a.e.\ on $\
\partial\mathbb D$\ .}\end{equation}

Moreover, by Remark 5 $\tilde f$ is a generalized analytic function
with a source $\tilde g\in L^p(\mathbb D)$,  $\tilde f={\cal
H}^*_{\tilde g}$, where ${\cal H}^*_{\tilde g}$ is the Hilbert
operator described in Section 3, Lemma 1, and associated with
$\tilde\lambda$ and $\tilde\varphi$, and the support of $\tilde g$
is in the support of $\tilde h$ and the upper bound of $\|\tilde
g\|_p$ depends only on $\|\tilde h\|_p$ and on the function $q$.

In addition, $\tilde g:\mathbb D\to\mathbb C$ is a fixed point of
the nonlinear operator $\tilde\Omega_{g_*}:=\tilde h\cdot q({\cal
H}^*_{g_*}):L^p_{\tilde h}(\mathbb D)\to L^p_{\tilde h}(\mathbb D)$,
where $L_{\tilde h}^{p}(\mathbb D)$ consists of functions $g_*$ in
$L^{p}(\mathbb D)$ with supports in the support of $\tilde h$.

Next, setting $f=\tilde f\circ c$, by simple calculations, see e.g.
Section 1.C in \cite{Alf}, we obtain that $\frac{\partial
f}{\partial\overline\xi}=\frac{\partial\tilde f}{\partial\overline
z}\circ c\cdot\overline{c^{\prime}}$ and, consequently, the function
$f:D\to\mathbb C$ is in the class $W^{1,p}_{\rm loc}\cap
C^{\alpha}_{\rm loc}(D)$ with $\alpha=(p-2)/p$ and satisfies
equation (\ref{eqQUASILINEARW5}). Moreover, $f$ is a generalized
analytic function with the source $g=\tilde g\circ c$ in the class
$L^p(D)$, $f(\xi)={\cal H}^*_{\tilde g}(c(\xi))$, and the support of
$g$ is in the support of $h$ and the upper bound of $\| g\|_p$
depends only on $\| h\|_p$, the function $q$ and the domain $D$.

It remains to show that $f$ has the angular limits as
$\xi\to\omega\in\partial D$ and satisfies the boundary condition
(\ref{eqLIMHh5}) a.e. on $\partial D$. Indeed, by the Lindel\"of
theorem, see e.g. Theorem II.C.2 in \cite{Ko}, if $\partial D$ has a
tangent at a point $\omega$, then $\arg\ [c_*(\omega)-c(\xi)]-\arg\
[\omega-\xi]\to\mathrm {const}$ as $\xi\to\omega$. In other words,
the images under the conformal mapping $c$ of sectors in $D$ with a
vertex at $\omega\in\partial D$ is asymptotically the same as
sectors in $\mathbb D$ with a vertex at
$\zeta=c_*(\omega)\in\partial \mathbb D$. Consequently,
nontangential paths in $D$ are transformed under $c$ into
nontangential paths in $\mathbb D$ and inversely a.e. on $\partial
D$ and $\partial{\mathbb D},$ respectively, because the rectifiable
boundary $\partial D$ has a tangent a.e. and $c_*$ and $c_*^{-1}$
keep sets of the length zero.
\end{proof}

\medskip

In particular, choosing $\lambda\equiv 1$ in Theorem 2 we obtain the
following consequence on the Dirichlet problem for the nonlinear
equations of the Vekua type.

\medskip

{\bf Corollary 2.} {\it Let  $D$ be  a Jordan domain with a
rectifiable boundary, $\varphi:\partial D\to\mathbb{R}$ be
measurable, $h: D\to\mathbb C$ be in $L^p(D)$, $p>2$, with compact
support in $D$, and let $q:\mathbb C\to\mathbb C$ be a continuous
function with condition (\ref{eqAPRIORYW5}).

Then there is $f: D\to\mathbb C$ in the class $W^{1,p}_{\rm loc}\cap
C^{\alpha}_{\rm loc}(D)$ with $\alpha=(p-2)/p$, satisfying equation
(\ref{eqQUASILINEARW5}), and the angular limits
\begin{equation}\label{eqLIMHhD} \lim\limits_{\xi\to\omega}\ \mathrm
{Re}\  f(\xi)\ =\ \varphi(\omega) \quad\quad\quad \mbox{a.e.\ on $\
\partial D$}\ .\end{equation}}

{\bf Remark 6.} Moreover, by the proof of Theorem 2 $f$ is a
generalized analytic function with a source $g\in L^p(D)$ whose
support is in the support of $h$ and the upper bound of $\| g\|_p$
depends only on $\| h\|_p$, the function $q$ and the domain $D$.

In addition, $g=\tilde g\circ c$ and $f={\cal H}^*_{\tilde g}\circ
c$, where $c$ is a conformal mapping of $D$ onto $\mathbb D$,
$\tilde g:\mathbb D\to\mathbb C$ is a fixed point of the nonlinear
operator $\tilde\Omega_{g_*}:=\tilde h\cdot q({\cal
H}^*_{g_*}):L^p_{\tilde h}(\mathbb D)\to L^p_{\tilde h}(\mathbb D)$,
where $L_{\tilde h}^{p}(\mathbb D)$ consists of functions $g_*$ in
$L^{p}(\mathbb D)$ with supports in the support of $\tilde h
:=h\circ C\cdot \overline{C^{\prime}}$, $C=c^{-1}$, ${\cal
H}^*_{\tilde g}$ is the Hilbert operator described in Section 3 and
associated with $\tilde\lambda=\lambda\circ c_*^{-1}$ and
$\tilde\varphi=\varphi\circ c_*^{-1}$. Here $c_*:\partial
D\to\partial\mathbb D$ is the homeomorphic boundary correspondence
under the mapping $c$.

\section{On completely continuous Poincare operators}

In Section 7 of \cite{R7}, we considered the Poincare boundary-value
problem on the directional derivatives and, in particular, the
Neumann problem with arbitrary measurable boundary data over the
natural parameter for the Poisson equations
\begin{equation}\label{eqSSC}
\triangle U(z)\ =\ G(z)
\end{equation}
with real valued functions $G$ of classes $L^{p}(D)$ with $p>2$ in
Jordan's domains $D$ in $\mathbb C$ with rectifiable boundaries.

Recall that a continuous solution $U$ of (\ref{eqSSC}) in the class
$W^{2,p}_{\rm loc}$ was called in \cite{R7} a {\bf generalized
harmonic function with the source} $\bf G$ and that by the Sobolev
embedding theorem such a solution belongs to the class $C^1$, see
Theorem I.10.2 in \cite{So}.

As usual, here $\frac{\partial u}{\partial \nu}\ (\xi)$  denotes the
derivative of $u$ at the point $\xi\in D$ in the direction
$\nu\in\mathbb C$, $|\nu|=1$, i.e.,
\begin{equation}\label{eqDERIVATIVE}
\frac{\partial u}{\partial \nu}(\xi)\ :=\ \lim_{t\to 0}\
\frac{u(\xi+t\cdot\nu)-u(\xi)}{t}\ .
\end{equation}
The  Neumann boundary value problem is a special case of the
Poincare problem on the directional derivatives with the unit
interior normal $n=n(\omega)$ to $\partial D$ at the point $\omega$
as $\nu(\omega)$, see Corollary 4 further.


By Theorem 6 in  \cite{R7}, for each measurable functions
$\nu:\partial\mathbb D\to\mathbb{C},\: |\nu(\zeta)|\equiv 1$ and
$\varphi:\partial\mathbb D\to\mathbb{R}$, $G:\mathbb D\to \mathbb R$
in $L^p(\mathbb D)$, $p>2$, there is a generalized harmonic function
$U:\mathbb D\to\mathbb R$ with the source $G$ that have the angular
limits
\begin{equation}\label{eqLIMITC} \lim\limits_{z\to\zeta}\ \frac{\partial U}{\partial \nu}\ (z)\ =\
\varphi(\zeta) \quad\quad\quad \mbox{a.e.\ on $\ \partial\mathbb
D$}\ .\end{equation} Furthermore, the space of such functions $U$
has the infinite dimension.


As it follows from constructions in the proofs of Theorem 1 and 6 in
\cite{R7}, see especially (2.6) there, one of such functions $U$ can
be presented as a sum of the logarithmic (Newtonian) potential
${\cal N}_G$ of the source $G$,
\begin{equation} \label{eqIPOTENTIAL} N_{G}(z)\  :=\
\frac{1}{2\pi}\int\limits_{\mathbb C} \ln|z-w|\, G(w)\ d\, m(w)\
,\end{equation} where $d\, m(w)$ corresponds to the Lebesgue measure
in the plane, i.e., the area, and the harmonic function
\begin{equation}\label{eqHARMONIC}\gamma(z) \ :=\ {\rm
Re}\,\int\limits_0^z\{\ {\cal H}^*_{G/2}(\xi)\ -\ T_{G/2}(\xi)\ \}\
d\,\xi\ ,\end{equation} where ${\cal H}^*_g$ is the Hilbert operator
described in Section 3 but with $\lambda=\bar\nu$ and where we
assumed that $G\in L^p(\mathbb D)$, $p>2$, with compact support in
$\mathbb D$.

Denoting by ${\cal P}^*_G$ the given correspondence between such
sources $G$ and the generalized harmonic functions with the sources
$G$ and the Poincare boundary condition (\ref{eqLIMITC}), we see
that ${\cal P}^*_G$ is a completely con\-ti\-nu\-ous operator over
each disk $|z|<r<1$ because the operators ${\cal H}^*_{G/2}$ and
$T_{G/2}$ are so and, in addition, the indefinite integral as well
as the operator of taking $\rm Re$ are bounded and linear. Thus, by
Lemma 1 and Remark 4 we come to the following statements.

\medskip

{\bf Lemma 2.}{\it\, Let $\nu:\partial\mathbb D\to\mathbb{C},\:
|\nu(\zeta)|\equiv1$, and $\varphi:\partial\mathbb D\to\mathbb{R}$
be measurable functions. Then there is a Poincare operator ${\cal
P}^*_G$ over the sources $G:\mathbb D\to \mathbb R$ in $L^p(\mathbb
D)$, $p>2$, with compact supports in $\mathbb D$, generating
generalized harmonic functions $U: \mathbb D\to\mathbb R$ with the
sources $G$ and the angular limits (\ref{eqLIMITC}), whose
restriction to sources $G$ with $\rm{supp}\,G\subseteq\mathbb
D_{\rho}$ is completely continuous over $\mathbb D_r$ for each
$\rho$ and $r\in(0,1).$ }

{\bf Remark 7.} Moreover, we may assume that the restriction of the
operator ${\cal P}^*_G$ to $\mathbb D_r$ under each $r\in(0,1)$ is
bounded at infinity in the sense that $\max\limits_{z\in\mathbb
D_r}|{\cal P}^*_G(z)|\leq M\cdot\|G\|_p$ for some $M>0$ and all $G$
with large enough $\| G\|_p$.

\section{On Poincare problem for semi-linear equations}

In this section we study the solvability of the Poincare
boundary-value problem for semi-linear Poisson equations of the form
$\triangle U(z) = H(z)\cdot Q(U(z))$ in the unit disk $\mathbb D$.
Again the Leray--Schauder approach allows us to reduce the problem
to the study of the linear Poisson equation from our last paper
\cite{R7} on the basis of Lemma 2 on completely continuous Poincare
operator ${\cal P}^*_G$ and Remark 7 on its boundedness at infinity
from the previous section.

\bigskip

Note that hypothesis (H2) in Section 2 will be automatically
satisfied in the proof of the next theorem because the initial
operator $F(\cdot):=F(\cdot,0)\equiv 0$ and hence $F$ has the only
one fixed point (at the origin) and its index is equal to $1$.

\bigskip

{\bf Theorem 3.} {\it Let $\nu:\partial\mathbb D\to\mathbb{C},\:
|\nu(\zeta)|\equiv1$, and $\varphi:\partial\mathbb D\to\mathbb{R}$
be measurable functions. Suppose that $H: \mathbb D\to\mathbb R$ is
a function in the class $L^p(\mathbb D)$ for $p>2$ with compact
support in $\mathbb D$ and $Q:\mathbb R\to\mathbb R$ is a continuous
function with
\begin{equation} \label{eqAPRIORYWH}
\lim\limits_{t\to \infty}\ \frac{Q(t)}{t}\ =\ 0\ .
\end{equation}

Then there is a function $U:\mathbb D\to\mathbb R$ in $W^{2,p}_{\rm
loc}(\mathbb D)\cap C^{1,\alpha}_{\rm loc}(\mathbb D)$ with
$\alpha=(p-2)/p$,
\begin{equation} \label{eqQUASILINEARWH}
\triangle U(z)\ =\ H(z)\cdot Q(U(z))\ \ \ \ \ \ \ \mbox{a.e. in
$\mathbb D$}\ ,
\end{equation}
and the angular limits
\begin{equation}\label{eqLIMHhH}
\lim\limits_{z\to\zeta}\ \frac{\partial U}{\partial \nu}\ (z)\ =\
\varphi(\zeta) \quad\quad\quad \mbox{a.e.\ on $\
\partial\mathbb D$}\ .\end{equation}}


\begin{proof} If $\| H\|_p=0$ or $\| Q\|_C= 0$, then any harmonic
function from Theorem 3 in \cite{R3} gives the desired solution of
(\ref{eqQUASILINEARWH}). Thus, we may assume that $\| H\|_p\neq 0$
and $\| Q\|_C\neq 0$. Set $Q_*(t)=\max\limits_{|\tau|\le
t}|Q(\tau)|$, $t\in\mathbb R^+:=[0,\infty)$. Then the function
$Q_*:\mathbb R^+\to\mathbb R^+$ is continuous and nondecreasing and,
moreover, by (\ref{eqAPRIORYWH})
\begin{equation} \label{eqAPRIORYWTH}
\lim\limits_{t\to \infty}\ \frac{Q_*(t)}{t}\ =\ 0\ .
\end{equation}

By Lemma 2 and Remark 7 we obtain the family of operators
$F(G;\tau): L_H^{p}(\mathbb D)\to L_H^{p}(\mathbb D)$, where
$L_H^{p}(\mathbb D)$ consists of functions $G\in L^{p}(\mathbb D)$
with supports in the support of $H$,
\begin{equation} \label{eqFORMULA}
F(G;\tau)\ :=\ \tau H\cdot Q({\cal P}^*_G) \ \ \ \ \ \ \ \forall\
\tau\in[0,1]
\end{equation}
which satisfies all groups of hypothesis H1-H3 of Theorem 1 in
\cite{LS}, see Proposition 1. Indeed:

H1). First of all, by Lemma 2 the function $F(G;\tau)\in
L_H^{p}(\mathbb D)$ for all $\tau\in[0,1]$ and  $G\in
L_H^{p}(\mathbb C)$ because the function $Q({\cal P}^*_G)$ is
continuous and, furthermore, the operators $F(\cdot\, ;\tau)$ are
completely continuous for each $\tau\in[0,1]$ and even uniformly
continuous with respect to the parameter $\tau\in[0,1]$.

H2). The index of the operator $F(\cdot\, ;0)$ is obviously equal to
$1$.

H3). Let us assume that solutions of the equations $G=F(G;\tau)$ is
not bounded in $L_H^{p}(\mathbb D)$, i.e., there is a sequence of
functions $G_n\in L_H^{p}(\mathbb D)$ with $\|G_n\|_p\to\infty$ as
$n\to\infty$ such that $G_n=F(G_n;\tau_n)$ for some
$\tau_n\in[0,1]$, $n=1,2,\ldots$. However, then by Remark 7 we have
that, for some constant $M>0$,
$$
\| G_n\|_p\ \le\ \| H\|_p\ Q_*\left(\, M\, \| G_n\|_p\right)\
$$
and, consequently,
\begin{equation} \label{eqEST3W}
\frac{Q_*(\, M\, \| G_n\|_p)}{M\, \| G_n\|_p}\ \ge\ \frac{1}{M\, \|
H\|_p}\ >\ 0
\end{equation}
for all large enough $n$. The latter is impossible by condition
(\ref{eqAPRIORYWTH}). The obtained contradiction disproves the above
assumption.

Thus, by Theorem 1 in \cite{LS} there is a function $G\in L_H^p(D)$
with $F(G;1)=G$, and by Lemma 2 the function $U:={\cal P}^*_G$ gives
the desired solution of (\ref{eqQUASILINEARWH}).
\end{proof}

{\bf Remark 8.} Moreover, by the proof of Theorem 3 the function $U$
is a generalized analytic function with a source $G\in L^p(\mathbb
D)$, $U={\cal P}^*_G$, where ${\cal P}^*_G$ is the Poincare operator
described in the last section, Lemma 2, and the support of $G$ is in
the support of $H$ and the upper bound of $\| G\|_p$ depends only on
$\| H\|_p$ and on the function $Q$.

In addition, the source $G:\mathbb D\to\mathbb C$ is a fixed point
of the nonlinear operator $\Omega_G:=h\cdot Q({\cal
P}^*_G):L^p_H(\mathbb D)\to L^p_H(\mathbb D)$, where
$L_H^{p}(\mathbb D)$ consists of functions $G$ in $L^{p}(\mathbb D)$
with supports in the support of $H$.

We are able to say more in Theorem 3 for the case of $\mathrm {Re}\
n(\zeta)\overline{\nu(\zeta)}>0$, where $n(\zeta)$ is  the inner
normal to $\partial\mathbb D$ at the point $\zeta$. Indeed, the
latter magnitude is a scalar product of $n=n(\zeta)$ and $\nu
=\nu(\zeta)$ interpreted as vectors in $\mathbb R^2$ and it has the
geometric sense of projection of the vector $\nu$ into $n$. In view
of (\ref{eqLIMHhH}), since the limit $\varphi(\zeta)$ is finite,
there is a finite limit $U(\zeta)$ of $U(z)$ as $z\to\zeta$ in
$\mathbb D$ along the straight line passing through the point
$\zeta$ and being parallel to the vector $\nu$ because along this
line
\begin{equation}\label{eqDIFFERENCEH} U(z)\ =\ U(z_0)\ -\ \int\limits_{0}\limits^{1}\
\frac{\partial U}{\partial \nu}\ (z_0+\tau (z-z_0))\ d\tau\
.\end{equation} Thus, at each point with condition (\ref{eqLIMHhH}),
there is the directional derivative
\begin{equation}\label{eqPOSITIVEH}
\frac{\partial U}{\partial \nu}\ (\zeta)\ :=\ \lim_{t\to 0}\
\frac{U(\zeta+t\cdot\nu)-U(\zeta)}{t}\ =\ \varphi(\zeta)\ .
\end{equation}

\bigskip

In particular, in the case of the Neumann problem, $\mathrm {Re}\
n(\zeta)\overline{\nu(\zeta)}\equiv 1>0$, where $n=n(\zeta)$ denotes
the unit inner normal to $\partial\mathbb D$ at the point $\zeta$,
and we have by Theorem 3 and Remark 8 the following significant
result.

\medskip

{\bf Corollary 3.} {\it Let $\varphi:\partial\mathbb D\to\mathbb{R}$
be measurable, $H:\mathbb D\to \mathbb R$ be in $L^p(\mathbb D)$,
$p>2$, with compact support in $\mathbb D$ and let $Q:\mathbb
R\to\mathbb R$ be a continuous function with condition
(\ref{eqAPRIORYWH}).

Then one can find generalized harmonic functions $U:\mathbb
D\to\mathbb R$ with a source $G\in L^p(\mathbb D)$ satisfying
equation (\ref{eqQUASILINEARWH}) such that a.e. on $\partial\mathbb
D$ there exist:

\bigskip

1) the finite limit along the normal $n(\zeta)$
$$
U(\zeta)\ :=\ \lim\limits_{z\to\zeta}\ U(z)\ ,$$

2) the normal derivative
$$
\frac{\partial U}{\partial n}\, (\zeta)\ :=\ \lim_{t\to 0}\
\frac{U(\zeta+t\cdot n(\zeta))-U(\zeta)}{t}\ =\ \varphi(\zeta)\ ,
$$

3) the angular limit
$$ \lim_{z\to\zeta}\ \frac{\partial U}{\partial n}\, (z)\ =\
\frac{\partial U}{\partial n}\, (\zeta)\ .$$ }

\section{On Poincare problem in Jordan domains with rectifiable boundaries}

{\bf Theorem 4.}{\it\, Let $D$ be  a Jordan domain with a
rectifiable boundary and let $\nu:\partial D\to\mathbb{C},\:
|\nu|\equiv 1$, and $\varphi:\partial D\to\mathbb{R}$ be measurable
over the natural parameter.

Suppose that $H: D\to\mathbb R$ is in $L^p(D)$ for $p>2$ with
compact support in $D$ and $Q:\mathbb R\to\mathbb R$ is a continuous
function with
\begin{equation} \label{eqAPRIORYW68}
\lim\limits_{t\to \infty}\ \frac{Q(t)}{t}\ =\ 0\ .
\end{equation}

Then there is a function $U:D\to\mathbb R$ in $W^{2,p}_{\rm
loc}(D)\cap C^{1,\alpha}_{\rm loc}(D)$ with $\alpha=(p-2)/p$,
\begin{equation} \label{eqQUASILINEARW568}
\triangle U(\xi)\ =\ H(\xi)\cdot Q(U(\xi))\ \ \ \ \ \ \ \mbox{a.e.
in $D$}\ ,
\end{equation}
and the angular limits
\begin{equation}\label{eqLIMIT68} \lim\limits_{\xi\to\omega}\ \frac{\partial U}{\partial \nu}\ (\xi)\ =\
\varphi(\omega) \quad\quad\quad \mbox{a.e.\ on $\ \partial D$}\
.\end{equation}}

\medskip

\begin{proof}
Arguing similarly to the first and second items in the proof of
Theorem 2, we see that $\tilde\nu :=\nu\circ c_*^{-1}$ and
$\tilde\varphi :=\varphi\circ c_*^{-1}$ are measurable over natural
parameters, where $c_*:=\tilde{c}|_{\partial D}: {\partial
D}\to\partial\mathbb D$ is the restriction to the boundary of the
homeomorphic extension $\tilde c$ of $c$ to $\overline{D}$ onto
$\overline{\mathbb D}$.

Now, set $\tilde H=|C^{\prime}|^2\cdot H\circ C$, where $C$ is the
inverse conformal mapping $C:=c^{-1}:\mathbb D\to D$. Then it is
clear by the hypotheses of Theorem 4 that $\tilde H$ has compact
support in $\mathbb D$ and belongs to the class $L^p(\mathbb D)$.
Consequently, by Theorem 3 there is $\tilde U:\mathbb D\to\mathbb R$
in $W^{1,p}_{\rm loc}(\mathbb D)\cap C^{1,\alpha}_{\rm loc}(\mathbb
D)$ with $\alpha=(p-2)/p$,
\begin{equation} \label{eqQUASILINEARW5568}
\triangle\tilde U(z)\ =\ \tilde H(z)\cdot Q(\tilde U(z))\ \ \ \ \ \
\ \mbox{a.e. in $\mathbb D$}
\end{equation}
and the angular limits
\begin{equation}\label{eqLIMHh5568}
\lim\limits_{z\to\zeta}\ \frac{\partial\tilde U}{\partial\tilde
\nu}\ (z)\ =\ \tilde\varphi(\zeta) \quad\quad\quad \mbox{a.e.\ on $\
\partial\mathbb D$\ .}\end{equation}

Moreover, by Remark 8 the function $\tilde U$ is a generalized
harmonic function with a source $\tilde G\in L^p(\mathbb D)$,
$\tilde U={\cal P}^*_{\tilde G}$, where ${\cal P}^*_{\tilde G}$ is
the Poincare operator described in Section 6, Lemma 2, and
associated with $\tilde\nu$ and $\tilde\varphi$, and the support of
$\tilde G$ is in the support of $\tilde H$ and the upper bound of
$\|\tilde G\|_p$ depends only on $\|\tilde H\|_p$ and on the
function $Q$.

In addition, $\tilde G:\mathbb D\to\mathbb C$ is a fixed point of
the nonlinear operator $\tilde\Omega_{G_*}:=\tilde H\cdot Q({\cal
P}^*_{G_*}):L^p_{\tilde H}(\mathbb D)\to L^p_{\tilde H}(\mathbb D)$,
where $L_{\tilde H}^{p}(\mathbb D)$ consists of functions $G_*$ in
$L^{p}(\mathbb D)$ with supports in the support of $\tilde H$.

Next, setting $U=\tilde U\circ c$, by simple calculations, see e.g.
Section 1.C in \cite{Alf}, we obtain that $\triangle
U=|c^{\prime}|^2\cdot\triangle\tilde U\circ c$  and, consequently,
the function $U:D\to\mathbb C$ is  in the class $W^{1,p}_{\rm
loc}(D)\cap C^{1,\alpha}_{\rm loc}(D)$ with $\alpha=(p-2)/p$ that
satisfies equation (\ref{eqQUASILINEARW568}), $U$ is a generalized
harmonic function with a source $G\in L^p(D)$ and, moreover,
$U(\xi)={\cal P}^*_{\tilde G}(c(\xi))$, where ${\cal P}^*_{\tilde
G}$ is the Poincare operator from Section 6, $G=\tilde G\circ c$,
and the support of $G$ is in the support of $H$ and the upper bound
of $\| G\|_p$ depends only on $\| H\|_p$, the function $Q$ and the
domain $D$.

Finally, arguing similarly to the last item in the proof of Theorem
2, we show that (\ref{eqLIMHh5568}) implies (\ref{eqLIMIT68}).
\end{proof}

{\bf Remark 9.} Moreover, by the proof of Theorem 4 the function $U$
is a generalized harmonic function with a source $G$ in the class
$L^p(D)$ whose support is in the support of $H$ and the upper bound
of $\| G\|_p$ depends only on $\| H\|_p$, the function $Q$ and the
domain $D$.

In addition, $G=\tilde G\circ c$ and $U={\cal P}^*_{\tilde G}\circ
c$, where $c$ is a conformal mapping of $D$ onto $\mathbb D$,
$\tilde G:\mathbb D\to\mathbb C$ is a fixed point of the nonlinear
operator $\tilde\Omega_{G_*}:=\tilde H\cdot Q({\cal
P}^*_{G_*}):L^p_{\tilde H}(\mathbb D)\to L^p_{\tilde H}(\mathbb D)$,
where $L_{\tilde H}^{p}(\mathbb D)$ consists of functions $G_*$ in
$L^{p}(\mathbb D)$ with supports in the support of $\tilde H:=H\circ
C\cdot \overline{C^{\prime}}$, $C=c^{-1}$, ${\cal P}^*_{\tilde G}$
is the Poincare operator described in Section 6 and associated with
$\tilde\nu=\nu\circ c_*^{-1}$ and $\tilde\varphi=\varphi\circ
c_*^{-1}$. Here $c_*:\partial D\to\partial\mathbb D$ is the
homeomorphic boundary correspondence under the mapping $c$.

We are able to say more in Theorem 4 for the case of $\mathrm {Re}\
n(\zeta)\overline{\nu(\zeta)}>0$, where $n(\zeta)$ is  the inner
normal to $\partial D$ at the point $\zeta$. Indeed, the latter
magnitude is a scalar product of $n=n(\zeta)$ and $\nu =\nu(\zeta)$
interpreted as vectors in $\mathbb R^2$ and it has the geometric
sense of projection of the vector $\nu$ into $n$. In view of
(\ref{eqLIMIT68}), since the limit $\varphi(\zeta)$ is finite, there
is a finite limit $U(\zeta)$ of $U(z)$ as $z\to\zeta$ in $D$ along
the straight line passing through the point $\zeta$ and being
parallel to the vector $\nu$ because along this line
\begin{equation}\label{eqDIFFERENCE} U(z)\ =\ U(z_0)\ -\ \int\limits_{0}\limits^{1}\
\frac{\partial U}{\partial \nu}\ (z_0+\tau (z-z_0))\ d\tau\
.\end{equation} Thus, at each point with condition
(\ref{eqLIMIT68}), there is the directional derivative
\begin{equation}\label{eqPOSITIVE}
\frac{\partial U}{\partial \nu}\ (\zeta)\ :=\ \lim_{t\to 0}\
\frac{U(\zeta+t\cdot\nu)-U(\zeta)}{t}\ =\ \varphi(\zeta)\ .
\end{equation}

\bigskip

In particular, in the case of the Neumann problem, $\mathrm {Re}\
n(\zeta)\overline{\nu(\zeta)}\equiv 1>0$, where $n=n(\zeta)$ denotes
the unit interior normal to $\partial D$ at the point $\zeta$, and
we have by Theorem 4 and Remark 9 the following significant result.

\medskip

{\bf Corollary 4.} {\it Let $D$ be a Jordan domain in $\Bbb C$ with
a rectifiable boundary and $\varphi:\partial D\to\mathbb{R}$ be
measurable over the natural parameter.

Suppose that $H: D\to \mathbb R$ is in $L^p(D)$, $p>2$, with compact
support in $D$. Then one can find a generalized harmonic function
$U: D\to\mathbb R$ with a source $G\in L^p(D)$ satisfying equation
(\ref{eqQUASILINEARW568}) such that a.e. on $\partial D$ there
exist:

\bigskip

1) the finite limit along $n(\zeta)$
$$
U(\zeta)\ :=\ \lim\limits_{z\to\zeta}\ U(z)\ ,$$

2) the derivative
$$
\frac{\partial U}{\partial n}\, (\zeta)\ :=\ \lim_{t\to 0}\
\frac{U(\zeta+t\cdot n(\zeta))-U(\zeta)}{t}\ =\ \varphi(\zeta)\ ,
$$

3) the angular limit
$$ \lim_{z\to\zeta}\ \frac{\partial U}{\partial n}\, (z)\ =\
\frac{\partial U}{\partial n}\, (\zeta)\ ,$$ where $n=n(\zeta)$ is
the unit inner normal at points $\zeta\in\partial D$.}

\section{The Poincare problem in physical applications}

Theorem 4 on the Poincare boundary-value problem with arbitrary
measurable boundary data over the natural parameter in Jordan
domains with rectifiable boundaries can be applied to mathematical
models of physical and chemical absorption with diffusion, plasma
states, stationary burning etc.

The first circle of such applications is relevant to
reaction-diffusion problems. Problems of this type are discussed in
\cite{Diaz}, p. 4, and, in detail, in \cite{Aris}. A nonlinear
system is obtained for the density $U$ and the temperature $T$ of
the reactant. Upon eliminating $T$ the system can be reduced to
equations of the type (\ref{eqQUASILINEARW568}),
\begin{equation}\label{RDP}
 \triangle  U\ =\ \sigma\cdot Q(U)
\end{equation}
with $\sigma >0$ and, for isothermal reactions, $Q(U) = U^{\beta}$
where $\beta>0$ that is called the order of the reaction. It turns
out that the density of the reactant $U$ may be zero in a subdomain
called a dead core. A particularization of results in Chapter 1 of
\cite{Diaz} shows that a dead core may exist just if and only if
$\beta\in(0,1)$ and $\sigma$ is large enough, see also the
corresponding examples in \cite{GNR}. In this connection, the
following statement may be of independent interest.

\medskip

{\bf Corollary 5.} {\it Let $D$ be  a Jordan domain in $\mathbb C$
with a rectifiable boundary, $\nu:\partial D\to\mathbb{C},\:
|\nu|\equiv 1$, and $\varphi:\partial D\to\mathbb{R}$ be measurable
over the natural parameter.

Suppose that $H: D\to\mathbb R$ is a function in the class $L^p(D)$
for $p>2$ with compact support in $D$.

Then there is a solution $U:D\to\mathbb R$  in the class
$W^{2,p}_{\rm loc}(D)\cap C^{1,\alpha}_{\rm loc}(D)$ with
$\alpha=(p-2)/p$ of the semi-linear Poisson equation
\begin{equation} \label{eqQUASILINEARW568C}
\triangle U(\xi)\ =\ H(\xi)\cdot U^{\beta}(\xi)\ ,\ \ \ 0\ <\ \beta\
<\ 1\  ,\ \ \ \ \ \ \ \mbox{a.e. in $D$}
\end{equation}
satisfying the Poincare boundary condition on directional
derivatives
\begin{equation}\label{eqLIMIT68C} \lim\limits_{\xi\to\omega}\ \frac{\partial U}{\partial \nu}\ (\xi)\ =\
\varphi(\omega) \quad\quad\quad \mbox{a.e.\ on $\ \partial
D$}\end{equation} in the sense of the angular limits. }

\medskip

Note also that certain mathematical models of a thermal evolution of
a heated plasma lead to nonlinear equations of the type (\ref{RDP}).
Indeed, it is known that some of them have the form
$\triangle\psi(u)=f(u)$ with $\psi'(0)=\infty$ and $\psi'(u)>0$ if
$u\not=0$ as, for instance, $\psi(u)=|u|^{q-1}u$ under $0 < q < 1$,
see e.g. \cite{Diaz}. With the replacement of the function
$U=\psi(u)=|u|^q\cdot {\rm sign}\, u$, we have that $u = |U|^Q\cdot
{\rm sign}\, U$, $Q=1/q$, and, with the choice $f(u) =
|u|^{q^2}\cdot {\rm sign}\, u$, we come to the equation $\triangle U
= |U|^q\cdot {\rm sign}\, U=\psi(U)$.

\medskip

{\bf Corollary 6.} {\it Let $D$ be  a Jordan domain in $\mathbb C$
with a rectifiable boundary, $\nu:\partial D\to\mathbb{C},\:
|\nu|\equiv 1$ and $\varphi:\partial D\to\mathbb{R}$ be measurable
over the natural parameter.

Suppose also that $H: D\to\mathbb R$ is a function in the class
$L^p(D)$ for $p>2$ with compact support in $D$.

Then there is a solution $U:D\to\mathbb R$  in the class
$W^{2,p}_{\rm loc}(D)\cap C^{1,\alpha}_{\rm loc}(D)$ with
$\alpha=(p-2)/p$ of the semi-linear Poisson equation
\begin{equation} \label{eqQUASILINEARW568CC}
\triangle U(\xi)\ =\ H(\xi)\cdot |U(\xi)|^{\beta -1} U(\xi)\ ,\ \ \
0\ <\ \beta\ <\ 1\  ,\ \ \ \ \ \ \ \mbox{a.e. in $D$}
\end{equation}
satisfying the Poincare boundary condition on directional
derivatives (\ref{eqLIMIT68C}).}

\bigskip

Finally, we recall that in the combustion theory, see e.g.
\cite{Barenblat}, \cite{Pokhozhaev} and the references therein, the
following model equation
\begin{equation}\label{combustion}
{\partial u(z,t)\over \partial t}\ =\ {1\over \delta}\cdot \triangle
u\ +\ e^{u}\ ,\ \ \ t\geq 0,\ z\in D,
\end{equation}
takes a special place. Here $u\ge 0$ is the temperature of the
medium and $\delta$ is a certain positive parameter. We restrict
ourselves here by the stationary case, although our approach makes
it possible to study the parabolic equation (\ref{combustion}), see
\cite{GNR}. Namely, the corresponding equation of the type
(\ref{eqQUASILINEARW568}) is appeared here after the replacement of
the function $u$ by $-u$ with the function $Q(u)=e^{-u}$ that is
bounded at all.

\medskip

{\bf Corollary 7.} {\it Let $D$ be  a Jordan domain in $\mathbb C$
with a rectifiable boundary, $\nu:\partial D\to\mathbb{C},\:
|\nu|\equiv 1$,  and $\varphi:\partial D\to\mathbb{R}$ be measurable
over the natural parameter.

Suppose also that $H: D\to\mathbb R$ is a function in the class
$L^p(D)$ for $p>2$ with compact support in $D$.

Then there is a solution $U:D\to\mathbb R$  in the class
$W^{2,p}_{\rm loc}(D)\cap C^{1,\alpha}_{\rm loc}(D)$ with
$\alpha=(p-2)/p$ of the semi-linear Poisson equation
\begin{equation} \label{eqQUASILINEARW568CCC}
\triangle U(\xi)\ =\ H(\xi)\cdot e^{U(\xi)}\ \ \ \ \ \ \ \mbox{a.e.
in $D$}
\end{equation}
satisfying the Poincare boundary condition on directional
derivatives (\ref{eqLIMIT68C}).}

\medskip

{\bf Remark 10.} In the above corollaries $U$ is a generalized
harmonic function with a source $G\in L^p(D)$ whose support is in
the support of $H$ and the upper bound of $\| G\|_p$ depends only on
$\| H\|_p$, the function $Q$ and the domain $D$.

In addition, $G=\tilde G\circ c$ and $U={\cal P}^*_{\tilde G}\circ
c$, where $c$ is a conformal mapping of $D$ onto $\mathbb D$,
$\tilde G:\mathbb D\to\mathbb C$ is a fixed point of the nonlinear
operator $\tilde\Omega_{G_*}:=\tilde H\cdot Q({\cal
P}^*_{G_*}):L^p_{\tilde H}(\mathbb D)\to L^p_{\tilde H}(\mathbb D)$,
where $L_{\tilde H}^{p}(\mathbb D)$ consists of functions $G_*$ in
$L^{p}(\mathbb D)$ with supports in the support of $\tilde H:=H\circ
C\cdot \overline{C^{\prime}}$, $C=c^{-1}$, ${\cal P}^*_{\tilde G}$
is the Poincare operator described in Section 6 and associated with
$\tilde\nu=\nu\circ c_*^{-1}$ and $\tilde\varphi=\varphi\circ
c_*^{-1}$. Here $c_*:\partial D\to\partial\mathbb D$ is the
homeomorphic boundary correspondence under the mapping $c$.

\section{Neumann problem in physical applications}

In turn, Corollary 4 can be applied to the study of the physical
phenomena discussed by us in the last section. In the connection,
the particular cases of the function $Q(t)$ of the forms $t^{\beta}$
and $|t|^{\beta -1}t$ with $\beta\in(0,1)$ and $e^t$ will be useful.

\medskip

{\bf Corollary 8.} {\it Let $D$ be a Jordan domain in $\Bbb C$ with
a rectifiable boundary and $\varphi:\partial D\to\mathbb{R}$ be
measurable over the natural parameter. Suppose that $H: D\to \mathbb
R$ is a function in the class $L^p(D)$, $p>2$, with compact support
in $D$.

Then one can find a generalized harmonic function $U: D\to\mathbb R$
with a source $G\in L^p(D)$ that is a solution $U:D\to\mathbb R$ in
the class $W^{2,p}_{\rm loc}(D)\cap C^{1,\alpha}_{\rm loc}(D)$ with
$\alpha=(p-2)/p$ of the semi-linear Poisson equation
\begin{equation} \label{eqQUASILINEARW568C1}
\triangle U(\xi)\ =\ H(\xi)\cdot U^{\beta}(\xi)\ ,\ \ \ 0\ <\ \beta\
<\ 1\  ,\ \ \ \ \ \ \ \mbox{a.e. in $D$}
\end{equation}
such that a.e. on $\partial D$ there exist:

1) the finite limit along the $n(\zeta)$
$$
U(\zeta)\ :=\ \lim\limits_{z\to\zeta}\ U(z)\ ,$$

2) the derivative
$$
\frac{\partial U}{\partial n}\, (\zeta)\ :=\ \lim_{t\to 0}\
\frac{U(\zeta+t\cdot n(\zeta))-U(\zeta)}{t}\ =\ \varphi(\zeta)\ ,
$$

3) the angular limit
$$ \lim_{z\to\zeta}\ \frac{\partial U}{\partial n}\, (z)\ =\
\frac{\partial U}{\partial n}\, (\zeta)\ ,$$ where $n=n(\zeta)$ is
the unit inner normal at points $\zeta\in\partial D$.}

\medskip

{\bf Corollary 9.} {\it Under hypotheses of Corollary 8, there is a
solution $U$ in the class $W^{2,p}_{\rm loc}(D)\cap
C^{1,\alpha}_{\rm loc}(D)$ with $\alpha=(p-2)/p$ of the semi-linear
Poisson equation
\begin{equation} \label{eqQUASILINEARW568CC1} \triangle U(\xi)\ =\
H(\xi)\cdot |U(\xi)|^{\beta -1} U(\xi)\ ,\ \ \ 0\ <\ \beta\ <\ 1\ ,\
\ \ \ \ \ \ \mbox{a.e. in $D$}
\end{equation}
such that a.e. on $\partial D$ all the conclusion 1)-3) of Corollary
8 hold, i.e., $U$ is a generalized solution of the Neumann problem
for (\ref{eqQUASILINEARW568CC1}) in the given sense.}

\medskip

{\bf Corollary 10.} {\it Under hypotheses of Corollary 8, there is a
solution $U$ in the class $W^{2,p}_{\rm loc}(D)\cap
C^{1,\alpha}_{\rm loc}(D)$ with $\alpha=(p-2)/p$ of the semi-linear
Poisson equation
\begin{equation} \label{eqQUASILINEARW568CCC1}
\triangle U(\xi)\ =\ H(\xi)\cdot e^{U(\xi)}\ \ \ \ \ \ \ \mbox{a.e.
in $D$}
\end{equation}
such that a.e. on $\partial D$ all the conclusion 1)-3) of Corollary
8 hold, i.e., $U$ is a generalized solution of the Neumann problem
for (\ref{eqQUASILINEARW568CCC1}) in the given sense.}




\noindent {\it Institute of Applied Mathematics and Mechanics \\ of
National Academy of Sciences of Ukraine, \\ Ukraine, Slavyansk},

\noindent Email: Ryazanov@nas.gov.ua

\bigskip

\noindent{\it Bogdan Khmelnytsky National University of Cherkasy,\\
Physics Dept., Lab. of Math. Phys.,\\
Ukraine, Cherkasy,}

\noindent Email: vl.ryazanov1@gmail.com

\end{document}